\documentclass{amsart}
\usepackage{amsmath,amsfonts,amsthm,amssymb,latexsym}
\usepackage{graphicx,float,wrapfig,epsf, caption, subcaption}
\usepackage{epsfig}
\usepackage{cyr}
\usepackage{psfrag}
\usepackage{tikz, xstring}
\usetikzlibrary{arrows,positioning,calc,fit,patterns,matrix,trees}

\newtheorem{thm}{Theorem}[section]
\newtheorem{cor}[thm]{Corollary}
\newtheorem{lem}[thm]{Lemma}

\newtheorem{prop}[thm]{Proposition}

\newtheorem{Example}[thm]{Example}

\def\<{\langle}
\def\>{\rangle}

\def\y{ {\text {\rm y}  } }

\def\F1{\mathcal{F}_1}

\def\.{\hskip.04cm}

\def\0{{\mathbf 0}}

\tikzset{
     nicearrow/.style={
     	->,
	thick,
	shorten <=2pt,
	shorten >=2pt,}
}

\makeatletter
\def\grd@save@target#1{%
  \def\grd@target{#1}}
\def\grd@save@start#1{%
  \def\grd@start{#1}}
\tikzset{
  grid with coordinates/.style={
    to path={%
      \pgfextra{%
        \edef\grd@@target{(\tikztotarget)}%
        \tikz@scan@one@point\grd@save@target\grd@@target\relax
        \edef\grd@@start{(\tikztostart)}%
        \tikz@scan@one@point\grd@save@start\grd@@start\relax
        \draw[major help lines] (\tikztostart) grid (\tikztotarget);
        \grd@start
        \pgfmathsetmacro{\grd@xa}{\the\pgf@x/1cm}
        \pgfmathsetmacro{\grd@ya}{\the\pgf@y/1cm}
        \grd@target
        \pgfmathsetmacro{\grd@xb}{\the\pgf@x/1cm}
        \pgfmathsetmacro{\grd@yb}{\the\pgf@y/1cm}
        \pgfmathsetmacro{\grd@xc}{\grd@xa + \pgfkeysvalueof{/tikz/grid with coordinates/major step}}
        \pgfmathsetmacro{\grd@yc}{\grd@ya + \pgfkeysvalueof{/tikz/grid with coordinates/major step}}
      }
    }
  },
  minor help lines/.style={
    help lines,
    step=\pgfkeysvalueof{/tikz/grid with coordinates/minor step}
  },
  major help lines/.style={
    help lines,
    line width=\pgfkeysvalueof{/tikz/grid with coordinates/major line width},
    step=\pgfkeysvalueof{/tikz/grid with coordinates/major step}
  },
  grid with coordinates/.cd,
  minor step/.initial=.2,
  major step/.initial=1,
  major line width/.initial=2pt,
}
\makeatother

\newcommand{\drawpermutate}[1]{%
    \begin{tikzpicture}
        \foreach [evaluate=\i as \x using int(\i-1)]\i in {0,1,...,8}
        {
            \foreach [evaluate=\j as \y using int(\j-1)] \j in {0,1,...,8}
            {
                \node at (\i,\j)[name=perm-\x-\y,]{};
            }
        }

        \draw ([xshift=-5mm]perm-1--1.south west)--([xshift=-5mm]perm-1-7.north west);
        \draw ([xshift=-5mm]perm-7--1.south west)--([xshift=-5mm]perm-7-7.north west);
        \draw ([yshift=5mm]perm--1-0.south west)--([yshift=5mm]perm-7-0.south east);
        \draw ([yshift=5mm]perm--1-6.south west)--([yshift=5mm]perm-7-6.south east);

        \foreach \mystyle/\coords/\leftorright in {#1}
        {
            \foreach \x/\y in \coords
            {
                \node[circle,fill=\mystyle,draw=\mystyle] at (perm-\x-\y){};
                \IfStrEq{\leftorright}{left}{%
                }
                {% otherwise put it on the right
                }
            }

        }

    \end{tikzpicture}
}

\newcommand{\shadetheboxes}[1]{
	\foreach \x/\y in {#1}
      	\fill[pattern color = black!65, pattern=north east lines] (\x,\y) rectangle +(1,1);
	}
	
\newcommand{\drawthegrid}[1]{
	\draw (0.01,0.01) grid (#1+0.99,#1+0.99);
	}
	
\newcommand{\drawtheclpattern}[1]{
	\foreach \x/\y in {#1}
      	\filldraw (\x,\y) circle (6pt);
	}
	
\newcommand{\drawtheclpatternwhite}[1]{
	\foreach \x/\y in {#1}
      	\draw[fill=white] (\x,\y) circle (6pt);
	}
	
\newcommand{\drawtheclpatternwhitebig}[1]{
	\foreach \x/\y in {#1}
		\draw[fill=white] (\x,\y) circle (11pt);
	}
	
\newcommand{\drawspecialbox}[1]{
	\foreach \x/\y/\z/\w/\A in {#1}
		{
       		\fill[color = white!100, opacity=1, rounded corners = 1.5pt] (\x+0.125,\y+0.125) rectangle (\z-0.125,\w-0.125);
       		\draw[color = black, rounded corners = 1.5pt] (\x+0.125,\y+0.125) rectangle (\z-0.125,\w-0.125);
       		\fill[black] (\x/2+\z/2,\y/2+\w/2) node {$\scriptstyle\A$};
       	}
    }

\newcommand{\drawspecialshadedbox}[1]{
	\foreach \x/\y/\z/\w/\A in {#1}
		{
       		\fill[color = white!100, opacity=1, rounded corners=1.5pt] (\x+0.125,\y+0.125) rectangle (\z-0.125,\w-0.125);
       		\fill[pattern color = black!65, pattern=north east lines, rounded corners=1.5pt] (\x+0.125,\y+0.125) rectangle (\z-0.125,\w-0.125);
       		\draw[color = black, rounded corners=1.5pt] (\x+0.125,\y+0.125) rectangle (\z-0.125,\w-0.125);
       		\fill[black] (\x/2+\z/2,\y/2+\w/2) node {$\scriptstyle\A$};
       	}
    }

\newcommand{\mpatternnl}[4]{										% mesh pattern - nl
  \raisebox{0.6ex}{
  \begin{tikzpicture}[scale=0.35, baseline=(current bounding box.center), #1]
  	\useasboundingbox (0.85,-0.1) rectangle (#2+1.4,#2+1.1);
	
    \shadetheboxes{#4}
    
    \drawtheclpattern{#3}
    
  \end{tikzpicture}}
}

\newcommand{\mmpattern}[5]{									% mesh patterns with a special box marked
  \raisebox{0.6ex}{
  \begin{tikzpicture}[scale=0.35, baseline=(current bounding box.center), #1]
  \useasboundingbox (0.0,-0.1) rectangle (#2+1.4,#2+1.1);
    
    \shadetheboxes{#4}
    
    \drawthegrid{#2}
    
    \drawspecialbox{#5}
    
    \drawtheclpattern{#3}

  \end{tikzpicture}}
}

\newcommand{\decpatternww}[8]{								% decorated pattern w white
  \raisebox{0.6ex}{
  \begin{tikzpicture}[scale=0.35, baseline=(current bounding box.center), #1]
  \useasboundingbox (0.0,-0.1) rectangle (#2+1.4,#2+1.1);
    
    \shadetheboxes{#6}
    
    \drawthegrid{#2}
       
    \drawspecialbox{#7}
    \drawspecialshadedbox{#8}
    
    \drawtheclpatternwhite{#4}
    \drawtheclpatternwhitebig{#5}
    \drawtheclpattern{#3}

  \end{tikzpicture}}
}

\pgfdeclareverticalshading{myshadingE}{80bp}
{color(-10bp)=(white); color(80bp)=(red)}

\begin{document}

\title[Enumeration of several two-by-four classes\\]
{Enumeration of several two-by-four classes}

\author[Sam~Miner]{ \ Sam~Miner$^\star$}

\thanks{\thinspace ${\hspace{-.45ex}}^\star$
Department of Mathematics,
Pomona College, Claremont, CA, 91711.
\hskip.06cm
Email:
\hskip.02cm
\texttt{samuel.miner@gmail.com}}

\date{\today}

\vskip1.3cm

\begin{abstract}
We use catalytic variables to derive generating functions for the permutation classes $Av(\textbf{4123},\textbf{1324}), Av(\textbf{4123},\textbf{1243})$, and $Av(\textbf{4123},\textbf{1342})$. Each generating function is algebraic of degree two, and the growth rates of the classes are 4, 5, and $\alpha$, respectively, where $\alpha \approx 4.17035$. As a consequence of our analysis, we see that a typical large permutation in $Av(\textbf{4123},\textbf{1324})$ is likely to also avoid \textbf{3124}, and is even more likely to avoid \textbf{31524}. Large permutations which avoid \textbf{4123} and \textbf{1243} are likely to also avoid \textbf{1423}.
\end{abstract}

\maketitle

\vskip.9cm

\section{Introduction}\label{s4:intro}
In this paper, we analyze permutation classes which avoid two patterns of length four, which we refer to as two-by-four classes. There are 56 two-by-four classes up to symmetry, and some of these symmetry classes are in fact Wilf-equivalent. Restricting our analysis to Wilf-equivalence classes means there are 38 classes to consider (see~\cite{Le}, \cite{Wik}). Of the Wilf-equivalence classes, 30 have been enumerated, and have been shown to have algebraic generating functions (see~\cite{AAB1}, \cite{AAB2}, \cite{AAV1}, \cite{AAV2}, \cite{Atk}, \cite{ASV}, \cite{Bev2}, \cite{Bev3}, \cite{Bo}, \cite{Cal}, \cite{K1}, \cite{K2}, \cite{KS}, \cite{Pan}, \cite{Vat1}). This leaves eight two-by-four classes which have yet to be enumerated. 

Albert et al.~\cite{AHPSV} argued that three of these eight unenumerated classes do not have differentially algebraic generating functions. The authors found functional equations that the generating functions must satisfy, and used these to enumerate the permutation classes for values of $n$ into the hundreds. They then showed that if the functions are differentially algebraic, they must be of a very large degree, or involve derivatives of a high degree, which is unlikely based on the generating functions which have been already been found for two-by-four classes.

Here, we give explicit generating functions for three of the five remaining two-by-four classes, namely $Av(\textbf{4123},\textbf{1324}), Av(\textbf{4123},\textbf{1243})$, and $Av(\textbf{4123},\textbf{1342})$. Each of these functions is algebraic of degree two. The growth rates of the classes are $4, 5$, and $\alpha$, respectively, where the growth rate $\alpha$ is approximately equal to $4.17035$. 

Our enumeration relies heavily on multivariate generating functions. By analyzing the structure of our classes, and by choosing catalytic variables appropriately, we derive equations that our generating functions must satisfy. Then, by making use of the kernel method, we solve these equations explicitly, finding closed forms for the desired functions. The structure of our argument mirrors closely one used by Bevan~\cite{Bev2} to enumerate the class $Av(\textbf{4213},\textbf{2413})$.

Our main theorems are given here.

\medskip

\begin{thm}\label{ss:1:t:1}
Let $\mathcal{P}_1(z)$ be the generating function for the permutation class $Av(\textbf{4123},\textbf{1324})$. Then we have
$$
\mathcal{P}_1(z) = \frac{2-13z+26z^2-17z^3+4z^4-z(1-2z-z^2)\sqrt{1-4z}}{2\sqrt{1-4z}(1-3z+z^2)^2}\,.
$$
\end{thm}

\medskip

\begin{thm}\label{s:2:t:2}
Let $\mathcal{P}_2(z)$ be the generating function for the permutation class $Av(\textbf{4123},\textbf{1243})$. Then we have
$$
\mathcal{P}_2(z) = \frac{2-8z+2z^2+17z^3-15z^4+4z^5-(2-10z+16z^2-9z^3+2z^4)\sqrt{1-6z+5z^2}}{2z(2-z)^2(1-2z)(1-3z+z^2)}\,.
$$ 
\end{thm}

\medskip

\begin{thm}\label{s:3:t:1}
Let $\mathcal{P}_3(z)$ be the generating function for the permutation class $Av(\textbf{4123},\textbf{1342})$. Then we have
$$
\mathcal{P}_3(z) = 1+ \frac{z(1-z)(1-2z)(1-7z+17z^2-16z^3+4z^4+(1-3z+3z^2)\sqrt{1-4z})}{2-22z+96z^2-220z^3+282z^4-196z^5+64z^6-8z^7}\,.
$$
\end{thm}

\medskip

In sections 2, 3, and 4, we discuss $Av(\textbf{4123},\textbf{1324}), Av(\textbf{4123},\textbf{1243})$, and $Av(\textbf{4123},\textbf{1342})$, respectively, and prove Theorems~\ref{ss:1:t:1},~\ref{s:2:t:2}, and~\ref{s:3:t:1}. We then make some final remarks and mention open problems, including the two remaining two-by-four classes.

For the remainder of the introduction, we introduce definitions and terminology we will need to explain our results. We assume some basic knowledge about pattern-avoidance and permutation classes, though the interested reader can see~\cite{Bev1} or~\cite{Vat2} for more a thorough discussion of the topic. 

We analyze permutations by considering their Hasse diagrams, and \emph{source graphs}, as described by Bevan~\cite{Bev3}. Each permutation can be decomposed into a series of source graphs in the following way: label the left-to-right minima in a permutation $\sigma$ by $m_1,m_2,\ldots,m_k$. Then, for each $i$, the source graph corresponding to $m_i$ consists of~$m_i$, and all elements above~$m_i$ which are not contained in $m_j$ for some $j<i$. These source graphs partition the vertices of~$\sigma$. 

We also use the term \emph{fan} to refer to a source graph which is a \textbf{123}-avoider. The term fan is appropriate, since if the elements are connected by edges if they form a \textbf{12}, as in the Hasse diagram, then every edge is incident to the left-to-right minimum.

As we also consider the asymptotic behavior of the number of permutations in these classes, we require the following proposition from analytic combinatorics. 

\begin{prop}\label{prop:asymp} (Flajolet \& Odlyzko \cite{FO}; see \cite{FS}, Theorem VI.1) 
The coefficient of $z^n$ in $\lambda(1-z\rho)^{\alpha}$ admits for large $n$ the following asymptotic expansion:
$$
[z^n]\lambda(1-z\rho)^{\alpha} \sim \frac{\lambda}{\Gamma(-\alpha)}\rho^n n^{-\alpha-1}\left(1+\sum_{k=1}^{\infty}\frac{e_k}{n^k}\right)\,,
$$
where
$$
e_k = \sum_{l=k}^{2k} \lambda_{k,l}\prod_{j=1}^l (\alpha+j) \qquad \text{ and }\qquad \lambda_{k,l} = [\nu^k
t^l] e^t(1+\nu t)^{-1-1/\nu}\,.
$$
\end{prop}

Recall that $[z^n]f(z)$ denotes the coefficient of $z^n$ in $f(z)$. Proposition~\ref{prop:asymp} enables us 
to determine the complete asymptotics of $[z^n]f(z)$ by repeated application to successive 
terms of the Puiseux expansion for f(z) around its dominant singularity.

\bigskip

\section{The permutation class $Av(\textbf{4123},\textbf{1324})$}\label{s:main}

\begin{figure}

\begin{tikzpicture}[scale=0.5]
\node at (0,0) {$Av(\textbf{4123}, \textbf{1324}, \textbf{1423}, \textbf{3124})$};
\draw[thick] (0.6,0.6) to (3.4,3.4);
\draw[thick] (-0.6,0.6) to (-3.4,3.4);
\draw[thick] (3.4,4.6) to (0.6,7.4);
\draw[thick] (-3.4,4.6) to (-0.6,7.4);
\draw[thick] (0,8.6) to (0,11.4);
\node at (-4,4) {$Av(\textbf{4123}, \textbf{1324}, \textbf{3124})$};
\node at (4,4) {$Av(\textbf{4123}, \textbf{1324}, \textbf{1423})$};
\node at (0,8) {$Av(\textbf{4123}, \textbf{1324}, \textbf{31524})$};
\node at (0,12) {$Av(\textbf{4123}, \textbf{1324})$};

\end{tikzpicture}
\caption{Subclasses within $Av(\textbf{4123}, \textbf{1324})$}
\label{fig:sets}
\end{figure}
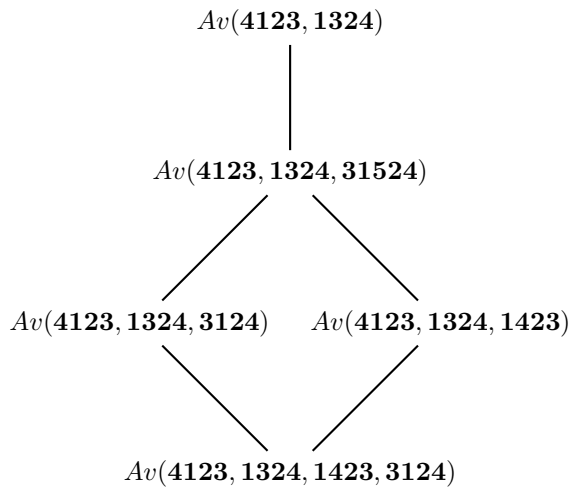

We first consider the permutation class $Av(\textbf{4123},\textbf{1324})$. We observe that there is a hierarchy of permutation classes within $Av(\textbf{4123},\textbf{1324})$, and enumerating the subclasses helps us enumerate the whole class. Concretely, 
$$
\aligned
&Av(\textbf{4123},\textbf{1324},\textbf{3124},\textbf{1423}) \subset Av(\textbf{4123},\textbf{1324},\textbf{1423}) \subset Av(\textbf{4123},\textbf{1324},\textbf{31524})\\
&\qquad \subset Av(\textbf{4123},\textbf{1324})\,.
\endaligned
$$ Additionally,  $$Av(\textbf{4123},\textbf{1324},\textbf{3124},\textbf{1423}) \subset Av(\textbf{4123},\textbf{1324},\textbf{3124}) \subset Av(\textbf{4123},\textbf{1324},\textbf{31524})\,.
$$ This subclass structure is shown in Figure~\ref{fig:sets}.

To enumerate these classes, we first state the following lemma. 

\begin{lem}\label{l:tgen}
Let $\mathcal{R}(t,z)$ be the generating function for a set of permutations, where $z$ measures the length, and $t$ the number of allowable positions for a new left-to-right minimum to be inserted. Let $\mathcal{C}(t,z)$ be the generating function for all permutations which can be generated by taking a permutation in $\mathcal{R}$, and adding a sequence of new source graphs which are fans, and whose non-minimum elements lie to the right of every previous element in the permutation. Then $\mathcal{C}(1,z) =\mathcal{R}(t,z)$, where $t = \frac{1-\sqrt{1-4z}}{2z}$ is the root of $1-t+t^2z=0$.
\end{lem}

This lemma helps us translate directly from bivariate generating functions to univariate generating functions, by substituting for $t$ with the appropriate function of $z$. We also remark that the function $t=\frac{1-\sqrt{1-4z}}{2z}$ is the generating function for the Catalan numbers, and for~$Av(\textbf{123})$. Bevan re-derives this well-known result in~\cite{Bev3}, and our argument is very similar to his.

\smallskip

\proof[Proof of Lemma~\ref{l:tgen}]
Every permutation counted by $\mathcal{C}$ is either in $\mathcal{R}$, or is formed by a permutation in $\mathcal{C}$ by adding a new fan as a source graph. This yields the functional equation $$\mathcal{C}(t,z) = \mathcal{R}(t,z) + \frac{tz}{1-tz}\frac{\mathcal{C}(1,z)-\mathcal{C}(t,z)}{1-t}\,.$$
The $\frac{tz}{1-tz}$ factor represents the new source graph which is added, and the $\frac{\mathcal{C}(1,z)-\mathcal{C}(t,z)}{1-t}$ factor represents the different allowable positions for the new minimum of the source graph to be placed.

Solving for $\mathcal{C}(t,z)$ gives us 
$$
((1-t)(1-tz)+tz)\mathcal{C}(t,z) = (1-t)(1-tz)\mathcal{R}(t,z) + tz\mathcal{C}(1,z)\,.
$$ 
Choosing $t=\frac{1-\sqrt{1-4z}}{2z}$ makes the left side equal to 0, and for this value of $t$, the right side simplifies to 
$$-tz\mathcal{R}(t,z)+tz\mathcal{C}(1,z)\,.
$$ 
Since the right side must also equal 0, we find that $\mathcal{C}(1,z) = \mathcal{R}(t,z)$, where 
$$
t = \frac{1-\sqrt{1-4z}}{2z}$$ is the generating function for the Catalan numbers, as desired. \qed 

\smallskip

Note that since $t$ satisfies $1-t+t^2z=0$, we also get $t^2=t-1$ and $\frac{1}{1-tz} = t$. We will use these alternate forms of expressing $t$ to simplify generating functions later on.

\medskip

\subsection{The permutation class $Av(\textbf{4123},\textbf{1324},\textbf{3124},\textbf{1423})$}\label{s:main:ss:4}

We first consider the class which avoids the patterns \textbf{4123}, \textbf{1324}, \textbf{3124}, and \textbf{1423}. Bruner analyzes this class in~\cite{Br}, and shows that $Av(\textbf{4123},\textbf{1324},\textbf{3124},\textbf{1423})$ is enumerated by the central binomial coefficients. Phrased in terms of the generating function for the class, we have the following theorem.

\smallskip

\begin{thm}\label{s:1:thm:br}{(Bruner, 2015)}
The generating function $\mathcal{A}(z)$ for $Av(\textbf{4123},\textbf{1324},\textbf{3124},\textbf{1423})$ has the form
$$
\mathcal{A}(z) = 1+\frac{z}{\sqrt{1-4z}}\,.
$$ 
\end{thm}

\smallskip

\proof
Here, we describe an alternate way to obtain the enumeration. We use $t$ to track the number of positions to the right of our current minimum where we are allowed to place a new minimum, as in~\cite{Bev2}. Our initial source graph has a leftmost element, plus an element of $Av(\textbf{213},\textbf{312})$ above it, since our permutation must avoid avoid $\textbf{1324}$ and $\textbf{1423}$. These initial source graphs are enumerated (see~\cite{SS}) by 

$$
1+ tz\left(1+\frac{tz}{1-2tz}\right) = 1+ \frac{tz(1-tz)}{1-2tz}\,.
$$ 

Now, when adding source graphs, avoiding \textbf{4123} means any subsequent source graph must be a fan. Since our permutations avoid \textbf{1423} and \textbf{1324}, there are no occurrences of \textbf{213} or \textbf{312} to the right of our current minimum element. Therefore, placing the new root of a source graph will never create a \textbf{3124} or \textbf{4123}, so the new root of a source graph can be placed anywhere to the right of our current minimum. Finally, since we avoid \textbf{3124}, any non-root elements of our source graphs must be to the right of everything currently in the permutation. 

Therefore, the hypotheses of Lemma~\ref{l:tgen} are satisfied, and we can substitute $t = \frac{1-\sqrt{1-4z}}{2z}$. We find that the generating function for the class $Av(\textbf{4123},\textbf{1324},\textbf{3124},\textbf{1423})$ is 
$$
\mathcal{A}(z) = 1+ \frac{tz(1-tz)}{1-2tz} = 1+ \frac{z}{1-2tz} = 1+ \frac{z}{\sqrt{1-4z}}\,,
$$
as desired.
\qed

\smallskip

These coefficients of this generating function are the central binomial coefficients, but even without that knowledge, we can calculate the asymptotic behavior of the number of permutations of length $n$ in this class, by using Proposition~\ref{prop:asymp}. 

\smallskip

\begin{cor}\label{s:1:cor:asymp1}
Asymptotically, $\vert \mathcal{S}_n \cap Av(\textbf{4123},\textbf{1324},\textbf{3124},\textbf{1423}) \vert \sim \frac{1}{4\sqrt{\pi}} 4^n n^{-\frac{1}{2}}$, as $n \to \infty$.
\end{cor}

\smallskip

\proof
Expanding $\mathcal{A}(z)$ around its dominant singularity, $\rho=\frac{1}{4}$, yields
$$
\mathcal{A}(z) \sim 1+\frac{1}{4}(1-4z)^{\frac{1}{2}} + O(1-4z)\,.
$$
Proposition~\ref{prop:asymp} then gives the desired result. 
\qed

\smallskip

We remark that for this class, we could also have derived the asymptotic behavior directly from Stirling's formula. On the other hand, Proposition~\ref{prop:asymp} is more broadly applicable, and is more effective when analyzing more complicated classes.

\begin{figure}
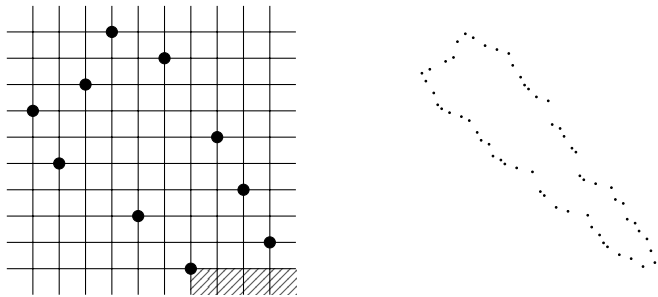

$$
\mmpattern{scale=1}{10}{1/7,2/5,3/8,4/10,5/3,6/9,7/1,8/6,9/4,10/2}{7/0,8/0,9/0,10/0}{}
\qquad\qquad
\mpatternnl{scale=.15}{60}{1/50, 2/48, 3/51, 4/45, 5/42, 6/41, 7/53, 8/40, 9/54, 10/58, 11/39, 12/60, 13/38, 14/59, 15/35, 16/33, 17/57, 18/32, 19/29, 20/56, 21/28, 22/27, 23/55, 24/52, 25/26, 26/49, 27/47, 28/46, 29/25, 30/44, 31/20, 32/19, 33/43, 34/37, 35/16, 36/36, 37/34, 38/15, 39/31, 40/30, 41/24, 42/23, 43/14, 44/11, 45/22, 46/9, 47/7, 48/6, 49/21, 50/18, 51/4, 52/17, 53/13, 54/3, 55/12, 56/10, 57/1, 58/8, 59/5, 60/2}{}
$$
\caption{Two permutations in $Av(\textbf{4123},\textbf{1324},\textbf{1423},\textbf{3124})$. The left permutation includes shaded positions where new left-to-right minima could be placed, while the right permutation is a larger example.}
\label{s:main:ss:4:fig:1}
\end{figure}

Two examples of permutations in $Av(\textbf{4123},\textbf{1324},\textbf{1423},\textbf{3124})$ are shown in Figure~\ref{s:main:ss:4:fig:1}\,.

\medskip

\subsection{The permutation class $Av(\textbf{4123},\textbf{1324},\textbf{3124})$}\label{s:main:ss:3}

We now consider a larger subclass, namely, permutations which can now contain occurrences of \textbf{1423}. The permutation class is $Av(\textbf{4123},\textbf{1324},\textbf{3124})$, and we have the following result.

\smallskip

\begin{thm}\label{ss:3:t:2}
Let $\mathcal{B}(z)$ be the generating function for $Av(\textbf{4123},\textbf{1324},\textbf{3124})$. Then 
$$
\mathcal{B}(z) = \frac{(1-3z)(1+\sqrt{1-4z})}{2\sqrt{1-4z}(1-3z+z^2)}\,.
$$
\end{thm}

\smallskip

\proof
We use the information from subsection~\ref{s:main:ss:4}, and enumerate the elements which are counted by $\mathcal{B}(z)$ but not $\mathcal{A}(z)$. These are permutations which \emph{contain} a \textbf{1423}. 

Let $\sigma$ be a permutation counted by $\mathcal{B}$ but not $\mathcal{A}$, so $\sigma$ contains an occurrence of \textbf{1423}. First, we argue that the leftmost element of $\sigma$, $\sigma_1$, must be the 1 in an occurrence of \textbf{1423}. Indeed, if $\sigma_1$ is above the 2 in the \textbf{1423} which occurs, then $\sigma_1$ either forms a \textbf{3124} or \textbf{4123} with the three smallest elements in the \textbf{1423}, both of which are forbidden. Therefore, $\sigma_1$ must be below the 2 in the \textbf{1423} which occurs, so must in fact serve as the 1 in a \textbf{1423}.

This observation implies that our initial source graph has a minimum, followed by an element of $Av(\textbf{4123},\textbf{213})$ which contains a \textbf{312}. Requiring our initial source graph to contain a \textbf{312} now restricts which positions are available for new minima. Specifically, anything left of the 1 but to the right of the 3 is disallowed for a new root of a source graph, to avoid creation of a~\textbf{4123}.

\begin{figure}
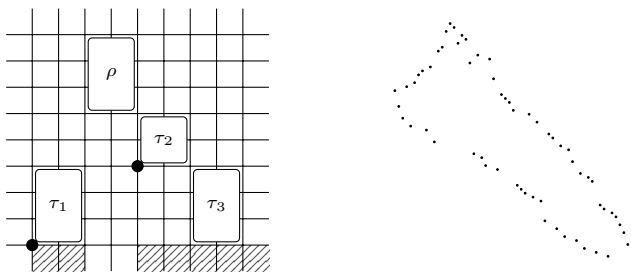

$$
\mmpattern{scale=1}{9}{1/1,5/4}{1/0,2/0,5/0,6/0,7/0,8/0,9/0}{1/1/3/4/\tau_1,3/6/5/9/\rho,5/4/7/6/\tau_2,7/1/9/4/\tau_3}
\qquad\qquad
\mpatternnl{scale=.15}{60}{1/43, 2/39, 3/36, 4/44, 5/34, 6/45, 7/47, 8/48, 9/33, 10/49, 11/30, 12/53, 13/54, 14/58, 15/60, 16/59, 17/55, 18/57, 19/56, 20/50, 21/27, 22/52, 23/26, 24/24, 25/51, 26/46, 27/23, 28/42, 29/41, 30/40, 31/38, 32/19, 33/18, 34/17, 35/37, 36/16, 37/35, 38/15, 39/10, 40/32, 41/31, 42/8, 43/29, 44/6, 45/28, 46/25, 47/5, 48/22, 49/3, 50/21, 51/20, 52/2, 53/14, 54/13, 55/1, 56/12, 57/11, 58/9, 59/7, 60/4}{}
$$
\caption{Two permutations in $Av(\textbf{4123},\textbf{1324},\textbf{3124})$. The left permutation has $\tau_1 \in Av(\textbf{21})$, $\tau_2,\tau_3 \in Av(\textbf{12}), \rho \in Av(\textbf{4123},\textbf{213})$, and $\rho$ and $\tau_2$ nonempty, and has the shaded cells representing potential positions of new left-to-right minima. The right is a random permutation of length 60.}
\label{s:1:fig:b}
\end{figure}

In Figure~\ref{s:1:fig:b} we display the typical structure of a permutation in $Av(\textbf{4123},\textbf{1324},\textbf{3124})$.

Since we want to track how many allowable insertion positions there are, we enumerate these initial source graphs by $$1 + tz \times \frac{1}{1-2tz} \times tz \times \frac{tz}{1-tz} \times \frac{z(1-z)}{1-3z+z^2}\,.$$

After the 1, our other terms come from the leftmost element, an increasing sequence merged with a decreasing sequence ($\tau_1$ and $\tau_3$), the rightmost 1 in a \textbf{312}, the elements which serve as twos corresponding to this 1, and a nonempty element of $Av(\textbf{4123},\textbf{213})$, each of which serves as a 3 in a \textbf{312}. 

Once we start adding new source graphs, we see that each subsequent source graph is again a fan, and all non-minimum elements must be to the right of the current permutation, since we are still avoiding \textbf{3124}. Therefore, we can apply Lemma~\ref{l:tgen}.

Applying Lemma~\ref{l:tgen} with the substitutions for $t$ yields the generating function 
$$\frac{t^4z^4(1-z)}{(1-3z+z^2)\sqrt{1-4z}}\,.$$ Therefore, combining this information with $\mathcal{A}(z)$ gives us 
$$
\mathcal{B}(z) = \mathcal{A}(z) + \frac{t^4z^4(1-z)}{(1-3z+z^2)\sqrt{1-4z}}\,,
$$
and simplifying gives us the result. \qed

\smallskip

As in Corollary~\ref{s:1:cor:asymp1}, Proposition~\ref{prop:asymp} helps us describe the asymptotic behavior of the size of our class.

\smallskip

\begin{cor}\label{s:1:cor:asymp2}
Asymptotically, $\vert \mathcal{S}_n \cap Av(\textbf{4123},\textbf{1324},\textbf{3124}) \vert \sim \frac{2}{5\sqrt{\pi}} 4^n n^{-\frac{1}{2}}$, as $n \to \infty$.
\end{cor}

\smallskip

\proof
Expanding $\mathcal{B}(z)$ around $\rho=\frac{1}{4}$, we get 
$$
\mathcal{B}(z) \sim \frac{2}{5}+\frac{2}{5}(1-4z)^{\frac{1}{2}}+O(1-4z)\,,
$$
and applying Proposition~\ref{prop:asymp} gives the result.
\qed

\smallskip

The first few terms of the series expansion of $\mathcal{B}(z)$ are 
$$
\mathcal{B}(z) = 1+z+2z^2+6z^3+21z^4+78z^5+297z^6+1143z^7+4419z^8+17119z^9+66836z^{10}+\ldots.
$$
Subsequent terms can be found at A277221 in the OEIS~\cite{Slo}.

\medskip

\subsection{The permutation class $Av(\textbf{4123},\textbf{1324},\textbf{31524})$}

Next, we consider the permutation class $Av(\textbf{4123},\textbf{1324},\textbf{31524})$. By avoiding \textbf{31524} while allowing \textbf{1423} and \textbf{3124}, we have a larger class of permutations. As it turns out, this subclass contains asymptotically almost every permutation that avoids \textbf{1423} and \textbf{1324}. We start with the following theorem.

\smallskip

\begin{thm}\label{ss:2:t:1}
Let $\mathcal{H}(z)$ be the generating function representing the permutation class $$Av(\textbf{4123},\textbf{1324},\textbf{31524})\,.$$ Then 
$$
\mathcal{H}(z) = \frac{3-22z+54z^2-54z^3+25z^4-4z^5-(1-6z+14z^2-16z^3+5z^4)\sqrt{1-4z}}{2\sqrt{1-4z}(1-3z+z^2)^2}\,.
$$
\end{thm}

\smallskip

\proof
We analyze the class by considering whether the permutations counted by $\mathcal{H}(z)$ additionally avoid \textbf{1423}. First, $Av(\textbf{4123},\textbf{1324},\textbf{31524},\textbf{1423}) = Av(\textbf{4123},\textbf{1324},\textbf{1423})$, since \textbf{31524} contains a \textbf{1423}. Also, since the basis elements in  $\{\textbf{4123},\textbf{1324},\textbf{1423}\}$ can be obtained from those in the set $\{\textbf{4123},\textbf{1324},\textbf{3124}\}$ by reversing, taking the inverse, and reversing each basis element again. Therefore, the subclass $Av(\textbf{4123},\textbf{1324},\textbf{1423})$ also has $\mathcal{B}(z)$ as its generating function. 

To calculate $\mathcal{H}(z)$, it suffices to count permutations which \emph{do} contain a \textbf{1423}. We use $\mathcal{H}_1(z)$ as the generating function which counts this permutations, so $\mathcal{H} = \mathcal{B} + \mathcal{H}_1$. 

We first count permutations in $\mathcal{H}_1$ which consist of one source graph, and we use $\mathcal{C}(t,z)$ to enumerate these. In the proof of Theorem~\ref{ss:3:t:2}, we saw that $$\mathcal{C}(t,z) = 1 + tz \times \frac{1}{1-2tz} \times tz \times \frac{tz}{1-tz} \times \frac{z(1-z)}{1-3z+z^2}\,.$$

We want to eventually apply Lemma~\ref{l:tgen}, though we need to realize that occurrences of \textbf{3124} are now allowed. Our permutation avoids \textbf{4123}, so any new source graph we add must still consist of a fan. Also, avoiding \textbf{1324} means that the third element of our fan must go to the right of every previous element, However, there is some flexibility on the positioning of our second element, though we can't have it form the 2 in an occurrence of \textbf{1423}. This basically means that we can add source graphs where the first non-minimum elements of our fans are directly adjacent to the minimum, or where they are to the right of every previous element. If a second element is neither, it then forms the 2 in a \textbf{1324} or a \textbf{1423}. 

Incorporating these new source graphs which are minima and elements directly next to them is allowed as long as the elements do not form a \textbf{4123}, meaning that they do not split any occurrence of \textbf{21} in the previous permutation. Also, since the previous permutation avoids \textbf{1324}, this means that reading left-to-right, once an element is the 2 in a \textbf{21}, the final element must be below this 2. In other words, the only allowable positions for placing a new minimum and adjacent element are to the left of the first occurrence of \textbf{21}. 

Therefore, we track how many elements are to the left of the first occurrence of \textbf{21}, with the variable $r$, and we track how many other elements do not form the 3 in a \textbf{312}, with the variable~$t$. Our initial source graphs have the generating function 
$$\mathcal{C}(r,t,z) = rz \times \left(\frac{1}{1-rz}\right) \times \frac{z-z^2}{1-3z+z^2} \times \frac{t^2z^2}{1-tz} \times \left(1+\frac{tz}{1-2tz}\right)\,.
$$
They have an initial element, an initial increasing sequence counted by powers of $r$, a nonempty member of $Av(\textbf{4123},\textbf{213})$ where each element serves as the 3 in a \textbf{312} (not counted by powers of $r$ or $u$), a \textbf{12} below and to the right of all these threes, with potentially multiple elements serving as the 2, and then potentially an interacting decreasing sequence below and to the right of the \textbf{312}, and an increasing sequence below and to the left of the elements we have just added. If these elements are here, the decreasing sequence must be nonempty, since otherwise our elements would be part of the initial increasing sequence.

Simplifying, we get 
$$\mathcal{C}(r,t,z) = \frac{rt^2z^4(1-z)}{(1-rz)(1-3z+z^2)(1-2tz)}\,.$$

After a sequence of source graphs completely contained among the elements counted by powers of $r$, our final source graph is either a single vertex, or a source graph of size 2, with a potentially empty sequence of new minima in between these elements. We use $\mathcal{D}(r,t,z)$ and~$\mathcal{E}(r,t,z)$ to denote these two separate cases. 

\begin{figure}
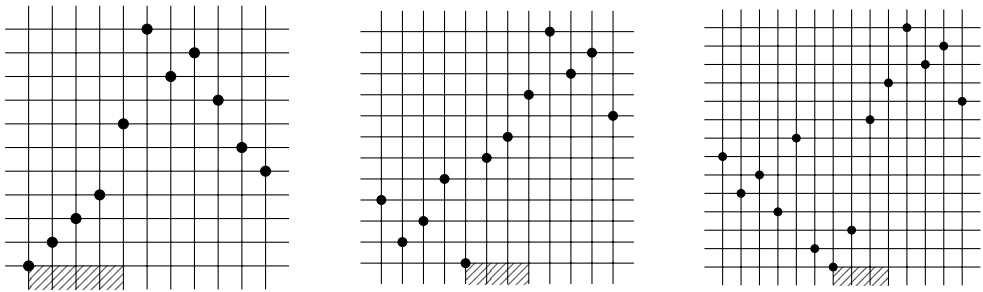

$$
\mmpattern{scale=0.9}{11}{1/1,2/2,3/3,4/4,11/5,10/6,5/7,9/8,6/11,7/9,8/10}{1/0,2/0,3/0,4/0}{}
\qquad
\mmpattern{scale=0.8}{12}{1/4, 2/2, 3/3, 4/5, 5/1, 6/6, 7/7, 8/9, 9/12, 10/10, 11/11, 12/8}{5/0,6/0,7/0}{}
\qquad
\mmpattern{scale=0.7}{14}{1/7, 2/5, 3/6, 4/4, 5/8, 6/2, 7/1, 8/3, 9/9, 10/11, 11/14, 12/12, 13/13, 14/10}{7/0,8/0,9/0}{}
$$
\caption{Permutations in $Av(\textbf{4123},\textbf{1324},\textbf{31524})$ enumerated by $\mathcal{C},\mathcal{D}$, and $\mathcal{E}$, respectively. Shaded cells represent possible positions for new minima which would keep the permutation in $\mathcal{D}$ and $\mathcal{E}$.}
\label{s:1:fig:CDE}
\end{figure}

In Figure~\ref{s:1:fig:CDE} we see three permutations, counted by $\mathcal{C}, \mathcal{D}$, and $\mathcal{E}$, respectively. The shaded areas represent possible positions for new minima to be placed, while keeping the permutation in $\mathcal{D}$, or $\mathcal{E}$.

The best way to view $\mathcal{D}$ and $\mathcal{E}$ is to think of them as inflations of \textbf{213}. Formally, elements in $\mathcal{D} \cup \mathcal{E}$ have the form $$
\textbf{213} [\sigma, \rho, \tau]\,,
$$ where $\sigma$ is nonempty in $Av(\textbf{4123},\textbf{132})$, and $\tau$ is in $Av(\textbf{4123},\textbf{213})$, so is enumerated by $\frac{\mathcal{C}}{z}$. We use $\frac{\mathcal{C}}{z}$ because $\sigma$ now plays the role of the minimum in the permutations enumerated by $\mathcal{C}$. The difference between $\mathcal{D}$ and $\mathcal{E}$ comes from $\rho$. For $\mathcal{D}$, $\rho$ is a single element, while for $\mathcal{E}$, $\rho$ is a fan of size at least two. 

These observations give us the generating functions $\mathcal{D}$ and $\mathcal{E}$, namely $$\mathcal{D}(r,u,z) = \mathcal{C}(r,u,z) \times \frac{1-z}{1-3z+z^2} \times z\,,$$ 
and
$$\mathcal{E}(r,u,z) = \mathcal{C}(r,u,z) \times \frac{z^2}{1-z} \times \frac{1-z}{1-3z+z^2} = \mathcal{C}(r,u,z) \times \frac{z^2}{1-3z+z^2}\,.$$ We consider $\mathcal{D}$ and $\mathcal{E}$ separately, since we need to track how many ways we can then take these initial permutations and create permutations which can be counted by Lemma~\ref{l:tgen}. 

A permutation in $Av(\textbf{4123},\textbf{1324},\textbf{31524})$ which contains a \textbf{1423}, either is counted by $\mathcal{C}, \mathcal{D}$, or $\mathcal{E}$, or can be formed by starting with a permutation counted by $\mathcal{C},\mathcal{D}$, or $\mathcal{E}$ and then adding a sequence of source graphs (fans) which preclude the permutation from being counted by $\mathcal{D}$ or~$\mathcal{E}$. We can count these extra permutations, if we initialize $\mathcal{R}$ correctly in Lemma~\ref{l:tgen}. 

One subset of permutations counted by $\mathcal{H}_1(z)$ is those which either consist of $\mathcal{C},\mathcal{D}$, or~$\mathcal{E}$, or which start with a permutation counted by $\mathcal{C},\mathcal{D}$, or~$\mathcal{E}$, and which then have a nonempty sequence of source graphs added, whose roots are all in positions corresponding to powers of $t$. If we refer to these permutations as being counted by $\mathcal{F}(t,z)$, then applying Lemma~\ref{l:tgen} with
$$
\mathcal{R}(t,z) = \left(\mathcal{C}(1,t,z)+\mathcal{D}(1,t,z) + \mathcal{E}(1,t,z)\right)
$$
gives us 
$$\mathcal{F}(t,z) = \frac{t^2z^4(1-z)^2}{(1-3z+z^2)^2(1-2tz)}\,.$$

The only permutations counted by $\mathcal{H}_1(z)$ which are not counted by $\mathcal{F}$ are those which are initially in $\mathcal{D} \cup \mathcal{E}$, and which no longer remain in $\mathcal{D} \cup \mathcal{E}$ after the addition of one or more non-minimal element to the most recent source graph, which is all the way to the right of the current permutation. We use $\mathcal{G}(t,z)$ to count this set, and separate $\mathcal{G}$ into $\mathcal{G}_{\mathcal{D}} + \mathcal{G}_{\mathcal{E}}$ based on where our initial permutation is. 

In terms of generating functions, from $\mathcal{D}$, we get $$
\mathcal{G}_{\mathcal{D}}(t,z) = \mathcal{D}(t,t,z) \times \frac{tz}{1-tz}\,.$$
From $\mathcal{E}$, we get $$
\mathcal{G}_{\mathcal{E}}(t,z) = 
\mathcal{E}(t,t,z) \times \frac{t^2z}{1-tz}\,.$$

Note that we get a second power of $t$ in $\mathcal{E}$. This is because $\mathcal{E}$ ends with a source graph of size two normally. We want to include the position between our minimum and second element as an allowable position for the root of a new source graph.

Using $\mathcal{G}_{\mathcal{D}} + \mathcal{G}_{\mathcal{E}}$ as $\mathcal{R}$ in Lemma~\ref{l:tgen}, we obtain $$
\mathcal{G}(t,z) = \mathcal{G}_{\mathcal{D}}+\mathcal{G}_{\mathcal{E}} = \frac{t^6z^6(1-z)}{(1-3z+z^2)^2(1-2tz)}\left(1-z+tz\right)\,,$$
where again $t$ is the Catalan generating function.

Finally, calculating 
$$
\mathcal{H}(z) = \mathcal{B}(z)+\mathcal{F}(t,z) + \mathcal{G}(t,z)\,,
$$
substituting with $$t = \frac{1-\sqrt{1-4z}}{2z}\,,$$ and simplifying, gives us the desired result.
\qed

\smallskip

By using Proposition~\ref{prop:asymp}, we can describe the asymptotic behavior of the number of permutations in $Av(\textbf{4123},\textbf{1324},\textbf{31524})$.

\smallskip

\begin{cor}\label{s:1:cor:asymp3}
Asymptotically, 
$$\left\vert \mathcal{S}_n \cap Av(\textbf{4123},\textbf{1324},\textbf{31524}) \right\vert \sim \frac{16}{25\sqrt{\pi}}4^n n^{-\frac{1}{2}}\left(1-\frac{105}{64n}+o\left(\frac{1}{n}\right)\right)\,,$$ as $n \to \infty$.
\end{cor}

\smallskip

\proof
The series expansion of $\mathcal{H}(z)$ around $z=\frac{1}{4}$ gives us 
$$
\mathcal{H}(z) = -\frac{16}{25\sqrt{1-4z}}-\frac{37}{50}+\frac{21}{10}\sqrt{1-4z}+O(1-4z)\,.
$$
By applying Proposition~\ref{prop:asymp} to these terms, we get the result.
\qed

\smallskip

The first few terms of the enumeration of $Av(\textbf{4123},\textbf{1324},\textbf{31524})$ are given by 
$$
\mathcal{H}(z) = 1+z+2z^2+6z^3+22z^4+86z^5+343z^6+1374z^7+5497z^8+21926z^9+87176z^{10}+\ldots
$$
Subsequent terms can be found at A277222 in the OEIS~\cite{Slo}.

\smallskip

An example of a permutation enumerated by $\mathcal{H}$ is shown in Figure~\ref{s:main:ss:2:fig:1}\,.

\begin{figure}
$$
\mpatternnl{scale=.15}{60}{1/41, 2/42, 3/40, 4/43, 5/38, 6/39, 7/44, 8/36, 9/34, 10/33, 11/37, 12/30, 13/26, 14/24, 15/45, 16/47, 17/49, 18/23, 19/52, 20/55, 21/56, 22/60, 23/57, 24/59, 25/58, 26/54, 27/53, 28/50, 29/51, 30/21, 31/48, 32/20, 33/18, 34/46, 35/35, 36/32, 37/31, 38/16, 39/29, 40/15, 41/14, 42/28, 43/13, 44/11, 45/27, 46/10, 47/25, 48/22, 49/8, 50/19, 51/6, 52/17, 53/4, 54/12, 55/2, 56/9, 57/7, 58/5, 59/1, 60/3}{}
$$
\caption{A permutation of length 60 in $Av(\textbf{4123}, \textbf{1324}, \textbf{31524})$}
\label{s:main:ss:2:fig:1}
\end{figure}

\medskip

\subsection{The permutation class $Av(\textbf{4123},\textbf{1324})$}

Finally, we are ready to enumerate the class $Av(\textbf{4123},\textbf{1324})$.

\proof[Proof of Theorem~\ref{ss:1:t:1}]
First, since we have $\mathcal{H}(z)$ already, we can restrict our attention to permutations in the class~$Av(\textbf{4123},\textbf{1324})$ which \emph{contain} an occurrence of $\textbf{31524}$. Suppose $\sigma$ is a such a permutation with an occurrence of \textbf{31524}. Observe that the elements which form the 4 must all be decreasing, since otherwise they would be part of occurrences of \textbf{4123}. We refer to these elements as $\rho_1$. Elements which form the 3 must be in $Av(\textbf{132},\textbf{4123})$, and elements which form the 5 must be in $Av(\textbf{4123,213})$. We refer to them as $\sigma_1$ and $\sigma_2$, respectively.

Observe that the elements left of the 1 are all vertically between the 2 and $\rho_1$, since having an element above $\rho_1$ forms a \textbf{4123}, and having one below the 2 forms a \textbf{1324}. 

The one possibility for additional elements which are not new minima is to the right of everything, below the 2. These elements must be decreasing, and we refer to them as $\rho_2$. Also, there can only be one element serving as the 2, since an increase would imply a \textbf{4123}, and a decrease would imply a \textbf{1324}.

\begin{figure}
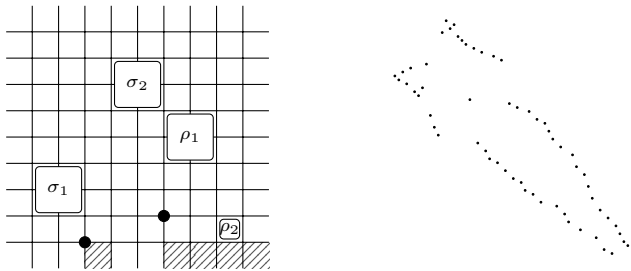

$$
\mmpattern{scale=1}{9}{3/1,6/2}{3/0,9/0,6/0,7/0,8/0}{1/2/3/4/\sigma_1,4/6/6/8/\sigma_2,6/4/8/6/\rho_1,8/1/9/2/\rho_2}
\qquad \qquad
\mpatternnl{scale=0.15}{60}{1/46, 2/45, 3/47, 4/44, 5/48, 6/42, 7/41, 8/43, 9/49, 10/36, 11/33, 12/31, 13/57, 14/60, 15/58, 16/59, 17/56, 18/55, 19/54, 20/40, 21/53, 22/29, 23/52, 24/27, 25/25, 26/51, 27/24, 28/50, 29/22, 30/39, 31/20, 32/38, 33/19, 34/17, 35/37, 36/16, 37/35, 38/14, 39/34, 40/32, 41/30, 42/13, 43/28, 44/11, 45/9, 46/26, 47/23, 48/7, 49/21, 50/18, 51/15, 52/5, 53/12, 54/2, 55/10, 56/1, 57/8, 58/6, 59/4, 60/3}{}
$$
\caption{Two permutations in $Av(\textbf{4123},\textbf{1324})$ which \emph{contain} a \textbf{31524}. The left permutation has $\rho_1,\rho_2 \in Av(\textbf{12})$, $\sigma_1 \in Av(\textbf{4123},\textbf{132})$, $\sigma_2 \in Av(\textbf{4123},\textbf{213})$, with $\rho_2$ as the only possible empty one. The shaded cells represent allowable positions for new left-to-right minima. The permutation on the right is of size 60, in $Av(\textbf{4123},\textbf{1324})$ and contains a \textbf{31524}.}
\label{s:1:fig:i}
\end{figure}

The structure of permutations counted by $\mathcal{I}(z)$ is displayed in Figure~\ref{s:1:fig:i}, as well as a sample permutation in $Av(\textbf{4123,1324})$ which contains a \textbf{31524}.

We generate these permutations using Lemma~\ref{l:tgen}. Our initial class is 
$$
\mathcal{R}(t,z) = \left(\frac{z-z^2}{1-3z+z^2}\right)^2\times \frac{tz}{1-tz} \times t^2z^2 \times \frac{1}{1-tz}\,.
$$ 

The first term comes from $\sigma_1$ and $\sigma_2$, the second from $\rho_1$, the third from our 2 element and our initial minimum, and the fourth from the possibly empty $\rho_2$.

Simplifying, and substituting with $t=\frac{1-\sqrt{1-4z}}{2z}$, our generating function is 
$$\mathcal{I}(z) = \frac{t^5z^5(1-z)^2}{(1-3z+z^2)^2} = \frac{(1-z)^2\left(1-5z+5z^2-(1-3z+z^2)\sqrt{1-4z}\right)}{2(1-3z+z^2)^2}\,.$$

Finally, calculating $\mathcal{P}_1(z) = \mathcal{H}(z)+\mathcal{I}(z)$ and simplifying gives us the desired expression for~$\mathcal{P}_1(z)$. \qed

\smallskip

By using Proposition~\ref{prop:asymp}, we can describe the asymptotic behavior of the number of permutations in $Av(\textbf{4123},\textbf{1324})$.

\smallskip

\begin{cor}\label{s:1:cor:asympP1}
Asymptotically, 
$$\left\vert \mathcal{S}_n \cap Av(\textbf{4123},\textbf{1324} \right\vert \sim \frac{16}{25\sqrt{\pi}}4^n n^{-\frac{1}{2}}\left(1-\frac{15}{16n}+o\left(\frac{1}{n}\right)\right)\,,
$$ as $n \to \infty$.
\end{cor}

\smallskip

\begin{cor}\label{s:1:cor:ratios}
Let $\sigma_n$ be a random permutation chosen uniformly from $\mathcal{S}_n \cap Av(\textbf{4123},\textbf{1324})$. Let $p_n$ be the probability that $\sigma_n$ avoids \textbf{31524}. Then 
$$
p_n \sim 1-\frac{45}{64n-60}+o\left(\frac{1}{n}\right) \quad \text{ as }\quad n \to \infty.
$$ Additionally, the probability that $\sigma_n$ avoids \textbf{3124} approaches $\frac{5}{8}$, and the probability that $\sigma_n$ avoids both \textbf{1423} and \textbf{3124} approaches $\frac{25}{64}$.
\end{cor}

\smallskip

\proof[Proof of Corollaries~\ref{s:1:cor:asympP1} and~\ref{s:1:cor:ratios}.]
The series expansion of $\mathcal{P}_1(z)$ around $z=\frac{1}{4}$ is
$$
\mathcal{P}_1(z) = -\frac{16}{25\sqrt{1-4z}}-\frac{14}{25}+\frac{12}{10}\sqrt{1-4z} + O(1-4z)\,.
$$

Applying Proposition~\ref{prop:asymp} to $\mathcal{P}_1(z)$ to the terms of this expansion, and comparing the results to 
Corollaries~\ref{s:1:cor:asymp1}, \ref{s:1:cor:asymp2}, and~\ref{s:1:cor:asymp3} completes the proof.
\qed

\smallskip

The first few terms of the enumeration of $Av(\textbf{4123},\textbf{1324})$ are given by 
$$
\mathcal{P}_1(z) = 1+z+2z^2+6z^3+22z^4+87z^5+352z^6+1428z^7+5768z^8+23156z^9+92416z^{10}+\ldots
$$
Subsequent terms can be found at A165532 in the OEIS~\cite{Slo}.

%$$
%\mmpattern{scale=1}{11}{1/1,2/2,3/3,4/4,11/5,10/6,5/7,9/8,6/11,7/9,8/10}{1/9,1/10,2/9,3/9,4/9,5/9,6/9,2/10,3/10,4/10,5/10,6/10}{1.5/4/4.5/5/r^3,9.5/6/11.5/7/t^2}
%$$
\bigskip

\section{The Permutation Class $Av(\textbf{4123},\textbf{1243})$}

We now turn our attention to the class $Av(\textbf{4123},\textbf{1243})$. As with our previous class, we first enumerate a large subclass, and then enumerate the rest of the permutations. 

\smallskip

\subsection{The permutation class $Av(\textbf{4123},\textbf{1243},\textbf{1423})$} 

We start by enumerating the subclass $Av(\textbf{4123},\textbf{1243},\textbf{1423})$. Our result is the following theorem.

\smallskip

\begin{thm}\label{s2:t1}
Let $\mathcal{J}(z)$ be the generating function for the class $Av(\textbf{4123},\textbf{1243},\textbf{1423})$. Then $$\mathcal{J}(z) = \frac{1+z-\sqrt{1-6z+5z^2}}{2(2z-z^2)}\,.$$
\end{thm}

\smallskip

We mention that several other three-by-four classes are Wilf-equivalent to this one (see~\cite{Bev2}, \cite{BB}, \cite{BHV}, \cite{CMS}). Before proving the theorem, we first give a lemma that describes the structure of our class. 

\smallskip

\begin{lem}\label{s:2:lem:grid}
The permutation class $Av(\textbf{4123},\textbf{1243},\textbf{1423})$ is the generalized grid class given by $$\mathbb{G} = \begin{tabular}{ | c |}
\hline
Av(21)\\
\hline
Av(123)\\
\hline
\end{tabular}\,.
$$

\end{lem}

\smallskip

\proof
First, if $\sigma \in \mathbb{G}$, then $\sigma$ can contain no occurrence of \textbf{4123, 1423}, or \textbf{1243}, since for each pattern, the largest two elements form a \textbf{21}, while the bottom three form a \textbf{123}. This shows that $\mathbb{G} \subset Av(\textbf{4123},\textbf{1243},\textbf{1423})$. 

To show the other containment, let $\sigma \in Av(\textbf{4123},\textbf{1243},\textbf{1423})$, and let $r$ be the largest integer which is the 1 in an occurrence of \textbf{21}. In other words, the elements from $r+1$ to $n$ avoid the pattern $\textbf{21}$. We claim that with these $n-r$ elements in the top cell, and the remaining $r$ elements in the bottom cell, $\sigma \in \mathbb{G}$. Clearly, the top cell is an element in $Av(21)$. Suppose for contradiction that the bottom cell contains an occurrence of \textbf{123}, with $\sigma_i, \sigma_j, \sigma_k$. These three elements along with $r+1$ must form a \textbf{1234}, since $r+1$ is the largest, and the other three possible patterns are forbidden. But then the elements $\sigma_i, \sigma_j, r+1$, and $r$ form a \textbf{1243}, a contradiction. Therefore $Av(\textbf{4123},\textbf{1243},\textbf{1423}) = \mathbb{G}$. 
\qed

\smallskip

\subsection{Proof of Theorem~\ref{s2:t1}}

With the help of Lemma~\ref{s:2:lem:grid}, we are ready to prove the theorem. 

\smallskip

\proof
Lemma~\ref{s:2:lem:grid} allows us to describe permutations in $Av(\textbf{4123},\textbf{1243},\textbf{1423})$ uniquely by a sequence of operations which are slightly different from adding source graphs. We call them \emph{source flags}. When adding a new source flag to a permutation, we add a new minimum, a decreasing sequence below and to the right of our previous permutation, which is mixed with a new increasing sequence, above and to the right of our previous permutation. 

\begin{figure}
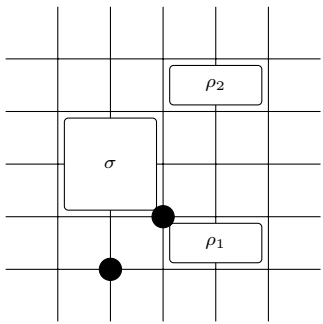


$$
\mmpattern{scale=2}{5}{2/1,3/2}{}{1/2/3/4/\sigma,3/4/5/5/\rho_2,3/1/5/2/\rho_1}
$$
\caption{Adding a new source flag to a permutation $\sigma$ in the class $Av(\textbf{4123},\textbf{1243},\textbf{1423})$. Here $\rho_1 \in Av(\textbf{12})$ and $\rho \in Av(\textbf{21})$, and their elements can be mixed horizontally, though the leftmost element must come from $\rho_1$.}
\label{s:2:fig:flag}

\end{figure}

We display the operation of adding a source flag in Figure~\ref{s:2:fig:flag}.

Observe that adding a new source flag to a permutation in $Av(\textbf{4123},\textbf{1243},\textbf{1423})$ will keep us in $Av(\textbf{4123},\textbf{1243},\textbf{1423})$. The new elements below the previous permutation form a fan, so the bottom of the new permutation is still a \textbf{123}-avoider. The new elements above the previous permutation still avoid \textbf{21}, so we remain in the grid class $\mathbb{G}$.

To enumerate permutations which can be created by adding source flags, we must track the number of allowable positions for new minima, which we do with a variable $v$. We consider the generating function 
$\mathcal{J}(v,z)$ where $z$ tracks the length of our permutation, and $v$ tracks the number of positions to the right of our current minimum.

Then $\mathcal{J}$ satisfies the functional equation 
$$
\mathcal{J}(v,z) = v+vz\left(1+\frac{vz}{1-2vz}\right)\left(\frac{\mathcal{J}(1,z)-\mathcal{J}(v,z)}{1-v}\right)\,.
$$

Solving for $\mathcal{J}(v,z)$, we get 
$$
(1-v-vz+2v^2z-v^2z^2)\mathcal{J}(v,z) = v(1-2vz)(1-v)+vz(1-vz)\mathcal{J}(1,z)\,.
$$
Choosing $v$ to cancel out the kernel on the left side, and solving for $\mathcal{J}(1,z)$, we find that
$$
\mathcal{J}(1,z) = v = \frac{1+z-\sqrt{1-6z+5z^2}}{2(2z-z^2)}\,,
$$
as desired. \qed

\smallskip

We can also use the generating function $\mathcal{J}(1,z)$ to compute the asymptotic behavior of the number of permutations in $Av(\textbf{4123},\textbf{1243},\textbf{1423})$ of length $n$. Bevan includes these calculations in Corollary 6 of~\cite{Bev2}, and we include the result here for reference.

\smallskip

\begin{cor}\label{s:2:c:asymJ}
Asymptotically, $\vert \mathcal{S}_n \cap Av(\textbf{4123},\textbf{1243},\textbf{1423}) \vert \sim \frac{5}{18}\sqrt{\frac{5}{\pi}} 5^n n^{-\frac{3}{2}}\,,$ as $n \to \infty$. 
\end{cor}

\smallskip

The first few terms of the series expansion of $\mathcal{J}(1,z)$ are
$$1+z+2z^2+6z^3+21z^4+79z^5+311z^6+1265z^7+5275z^8+22431z^9+96900z^{10}+\ldots
$$
Subsequent terms can be found at A033321 in the OEIS~\cite{Slo}.

\medskip

\subsection{The class $Av(\textbf{4123},\textbf{1243})$}
We are now ready to enumerate the whole class. 

\smallskip

\proof[Proof of Theorem~\ref{s:2:t:2}]
Using the information about $\mathcal{J}(1,z)$, it suffices to enumerate permutations which contain at least one occurrence of \textbf{1423}, but avoid \textbf{4123} and \textbf{1243}. We construct such permutations $\sigma$ in a systematic way. Let $\mathcal{K}(v,z)$ be the generating function for such permutations, where $n$ tracks length and $v$ tracks the number of allowable positions for a new left-to-right minimum to be placed. 

Permutations counted by $\mathcal{K}$ each contain at least one occurrence of \textbf{1423}. Consider the rightmost 3 in a \textbf{1423}, and the rightmost 2 corresponding to such a \textbf{1423}. Suppose there is just one such element which serves as the 1 in such a \textbf{1423}, and it is the smallest element of our permutation. Let $\mathcal{K}_1$ be the generating function for these permutations. We can describe permutations counted by $\mathcal{K}_1$ as consisting of a rightmost 3 in a \textbf{1423}, a lowest 2 in a \textbf{1423} involving the rightmost 3, the single minimum 1 in a \textbf{1423} involving this 2 and 3, and a nonempty permutation in $Av(132,1423)$ which plays the role of the 4, plus potentially more elements. There can only be one such 3, since a decrease where the 3 is would form a \textbf{1243}, and an increase where the 3 is would form a \textbf{4123}. There is room for a sequence of elements which are vertically between the 3 and the 2, and which are either left of our minimum, or vertically between our minimum and our 4. Also, there can be elements to the right of our 3, though they must either be increasing above everything, or decreasing between our 2 and our 1. 

\begin{figure}
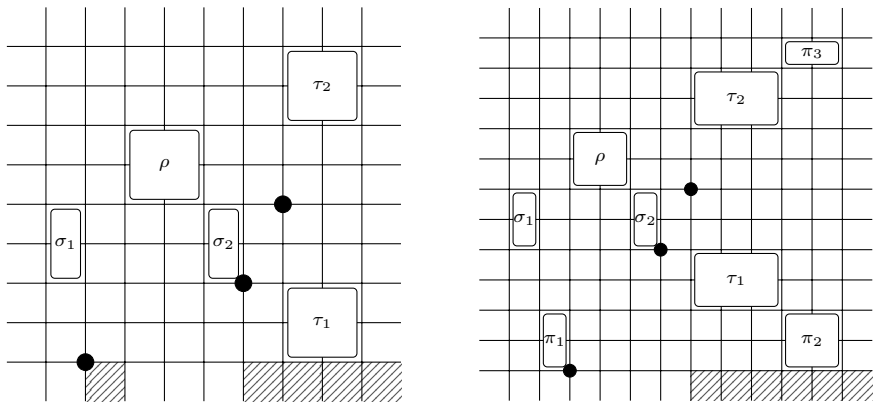


$$
\mmpattern{scale=1.5}{9}{2/1,6/3,7/5}{2/0,9/0,6/0,7/0,8/0}{1/3/2/5/\sigma_1,3/5/5/7/\rho,5/3/6/5/\sigma_2,7/1/9/3/\tau_1,7/7/9/9/\tau_2}
\qquad
\mmpattern{scale=1.15}{12}{3/1,7/7,6/5}{9/0,10/0,7/0,8/0,11/0,12/0}{1/5/2/7/\sigma_1,3/7/5/9/\rho,5/5/6/7/\sigma_2,7/3/10/5/\tau_1,7/9/10/11/\tau_2,2/1/3/3/\pi_1,10/1/12/3/\pi_2,10/11/12/12/\pi_3}
$$
\caption{The structure of two permutations counted by $\mathcal{K}_1$ and $\mathcal{K}_2$. Permutations $\sigma_1, \sigma_2, \tau_1,\pi_1,$ and $\pi_2$ are all in $Av(\textbf{12})$, all potentially empty. $\tau_2$ and $\pi_2$ are in $Av(\textbf{21})$, potentially empty, and $\rho \in Av(\textbf{132,4123})$, nonempty. The shaded cells represent the positions where new left-to-right minima could be inserted.}
\label{s:2:fig:k1}

\end{figure}

The structure of permutations counted by $\mathcal{K}_1$ is exhibited in Figure~\ref{s:2:fig:k1}. 

Combining all the elements of $\mathcal{K}_1(v,z)$, we get the generating function
$$
\mathcal{K}_{1}(v,z) = (vz)^3 \times \frac{z-z^2}{1-3z+z^2} \times \frac{1}{1-2z} \times \frac{1}{1-2vz} = \frac{v^3z^4(1-z)}{(1-3z+z^2)(1-2z)(1-2vz)}\,.
$$

We now can add new source flags to elements of $\mathcal{K}_1(v,z)$, to create new permutations counted by $\mathcal{K}$. 

First, consider a source flag whose minimum is directly next to the previous minimum, so still a 1 in our original rightmost \textbf{1423}. Adding a source flag of this type is equivalent to multiplying our generating function by a term of the form 
$$
z \left(1+\frac{vz}{1-2vz}\right)\,.
$$
Adding a sequence of these source flags results in the generating function 
$$
\mathcal{K}_2(v,z) = \frac{\mathcal{K}_1(v,z)}{v} \times \frac{1}{1-z\left(1+\frac{vz}{1-2vz}\right)} = \frac{v^2z^4(1-z)}{(1-3z+z^2)(1-2z)(1-(1+2v)z+vz^2)}\,.
$$

Observe that we divide by $v$ since we don't want to count the position directly next to our current minimum as an allowable position for a new minimum any more.

Finally, we can take permutations counted by $\mathcal{K}_2(v,z)$ and add source flags whose minima are not part of any occurrences of \textbf{1423}. Every permutation in $\mathcal{K}$ is either in $\mathcal{K}_2$, or it can be obtained from an element of $\mathcal{K}_2$ by a sequence of source flags. This gives us the functional equation

$$
\mathcal{K}(v,z) = \mathcal{K}_{2}(v,z)+vz\left(1+\frac{vz}{1-2vz}\right)\left(\frac{\mathcal{K}(1,z)-\mathcal{K}(v,z)}{1-v}\right)\,.
$$

Solving for $\mathcal{K}(1,z)$, and simplifying, we get 
$$
\mathcal{K}(1,z) = \frac{(1-z)(1-5z+4z^2-2z^3-(1-2z)\sqrt{1-6z+5z^2})}{2(1-2z)(1-3z+z^2)(2-z)^2}\,.
$$

Finally, calculating $\mathcal{P}_2(z)$ by adding $\mathcal{K}(1,z)$ and $\mathcal{J}(1,z)$ completes the proof of Theorem~\ref{s:2:t:2}. \qed

\smallskip

\begin{figure}
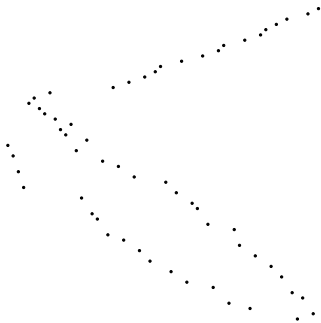

$$
\mpatternnl{scale=0.2}{60}{1/34, 2/32, 3/29, 4/26, 5/42, 6/43, 7/41, 8/40, 9/44, 10/39, 11/37, 12/36, 13/38, 14/33, 15/24, 16/35, 17/21, 18/20, 19/31, 20/17, 21/45, 22/30, 23/16, 24/46, 25/28, 26/14, 27/47, 28/12, 29/48, 30/49, 31/27, 32/10, 33/25, 34/50, 35/8, 36/23, 37/22, 38/51, 39/19, 40/7, 41/52, 42/53, 43/4, 44/18, 45/15, 46/54, 47/3, 48/13, 49/55, 50/56, 51/11, 52/57, 53/9, 54/58, 55/6, 56/1, 57/5, 58/59, 59/2, 60/60}{}
$$
\caption{The structure of permutations in $Av(\textbf{4123},\textbf{1243})$. The shading on the left permutation represents possible positions of new left-to-right minima. The right permutation is an example of size 60.}
\label{s:2:fig:2}
\end{figure}

An example of a permutation in $Av(\textbf{4123},\textbf{1243})$ of size 60 is displayed in Figure~\ref{s:2:fig:2}.

We can also use the generating function $\mathcal{P}_2(z)$ to compute the asymptotic behavior of the number of permutations in $Av(\textbf{4123},\textbf{1243})$ of length $n$.

\smallskip

\begin{cor}\label{s:2:c:asym}
Asymptotically, $$\vert \mathcal{S}_n \cap Av(\textbf{4123},\textbf{1243}) \vert \sim \frac{595}{1782}\sqrt{\frac{5}{\pi}} 5^n n^{-\frac{3}{2}}\,,$$ as $n \to \infty$. Therefore, the probability that a large permutation in $Av(\textbf{4123},\textbf{1243})$ also avoids \textbf{1423} tends to $\frac{99}{119}$.
\end{cor}

\smallskip

\proof
The proof follows from computing the Puiseux expansion for $\mathcal{P}_2(z)$ around the dominant singularity $\rho = \frac{1}{5}$, and applying Proposition~\ref{prop:asymp}. 

Then, using Corollary~\ref{s:2:c:asymJ}, we see that as $n$ gets large, the ratio of permutations in the class $Av(\textbf{4123},\textbf{1243})$ which also avoid \textbf{1423} behaves like
$$
\frac{\vert \mathcal{S}_n \cap Av(\textbf{4123},\textbf{1243},\textbf{1423}) \vert}{\vert \mathcal{S}_n \cap Av(\textbf{4123},\textbf{1243})\vert} \sim \frac{\frac{5}{18}\sqrt{\frac{5}{\pi}} 5^n n^{-\frac{3}{2}}}{\frac{595}{1782}\sqrt{\frac{5}{\pi}} 5^n n^{-\frac{3}{2}}} = \frac{(1782)(5)}{(595)(18)} = \frac{99}{119}\,,
$$
as desired.
\qed

\smallskip

In conclusion, most permutations in $Av(\textbf{4123},\textbf{1243})$ also avoid \textbf{1423}, though there is a nonzero proportion which contain \textbf{1423}. 

The first few terms of the series expansion of $\mathcal{P}_2(z)$ are
$$1+z+2z^2+6z^3+22z^4+88z^5+365z^6+1540z^7+6568z^8+28269z^9+122752z^{10}+\ldots
$$
Subsequent terms can be found at A165536 in the OEIS~\cite{Slo}.

\bigskip

\section{The Permutation Class $Av(\textbf{4123},\textbf{1342})$}

\subsection{Main result and structure of $Av(\textbf{4123},\textbf{1342})$.}

\begin{figure}
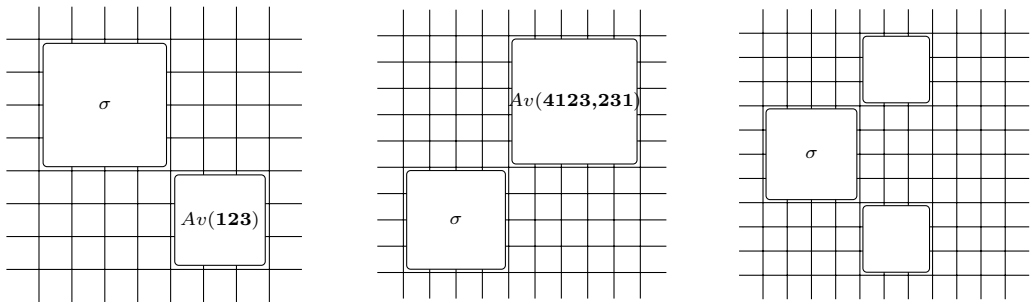

$$
\mmpattern{scale=1.25}{8}{}{}{1/4/5/8/\sigma,5/1/8/4/Av(\textbf{123})}
\qquad
\mmpattern{scale=1}{10}{}{}{1/1/5/5/\sigma,5/5/10/10/Av(\textbf{4123,231})}
\qquad
\mmpattern{scale=0.92}{11}{}{}{1/4/5/8/\sigma,5/1/8/4/,5/8/8/11/}
$$
\caption{The three ways of adding a new component to a permutation $\sigma$ in $Av(\textbf{4123},\textbf{1342})$. The diagrams correspond to $\mathcal{L},\mathcal{M}$, and $\mathcal{N}$, respectively.}
\label{s:3:fig:1}
\end{figure}

We now turn our attention to a third two-by-four class, the permutation class $Av(\textbf{4123},\textbf{1342})$. Before proving Theorem~\ref{s:3:t:1}, we first consider the structure of permutations in our class. We can construct permutations in this class in a unique way by a sequence of operations which are either: add a new maximal sum-component, add a new maximal skew-component, or add elements on the right which form the mix of a skew-component and a sum-component. This is exhibited in Figure~\ref{s:3:fig:1}. 

Let $\mathcal{L}(z)$ be the generating function for new nonempty skew-components, let $\mathcal{M}(z)$ be the generating function for new nonempty sum-components, and let $\mathcal{N}(z)$ be the generating functions for new nonempty minimal components which have elements both above and below the current permutation. By minimal, we mean that these components cannot be decomposed further into skew-components, sum-components, or smaller elements counted by $\mathcal{N}$. We refer to components counted by $\mathcal{N}$ as \emph{mixes}. 

By our requirement that elements counted by $\mathcal{L}$ and $\mathcal{M}$ be maximal, we get the following lemma.

\smallskip

\begin{lem}\label{s:3:lem:struct}
The generating function $\mathcal{P}_3(z)$ has the form 
$$
\mathcal{P}_3(z) = 1+z\frac{\left(\frac{(1+\mathcal{L})(1+\mathcal{M})}{1-\mathcal{L}\mathcal{M}}\right)}{1-\mathcal{N}\left(\frac{(1+\mathcal{L})(1+\mathcal{M})}{1-\mathcal{L}\mathcal{M}}\right)}\,.
$$
\end{lem}

\smallskip

\proof
We can describe the nonempty permutations in $\mathcal{P}_3(z)$ by a sequence of operations which either correspond to $\mathcal{L},\mathcal{M}$, or $\mathcal{N}$. Since $\mathcal{L}$ and $\mathcal{M}$ correspond to maximal skew- and sum-components, our sequence cannot have two consecutive $\mathcal{L}$ operations or two consecutive $\mathcal{M}$ operations. In other words, each nonempty permutation in $Av(\textbf{4123,1342})$ corresponds to a finite word with letters in $\{\mathcal{L},\mathcal{M},\mathcal{N}\}$, with no consecutive occurrences of $\mathcal{L}$ or $\mathcal{M}$, where each letter in fact represents a monomial corresponding to a permutation from its class.

Therefore, we can describe our sequence of operations as a potentially empty alternating sequence from $\{\mathcal{L},\mathcal{M}\}$, followed by a sequence of expressions which look like a mix (counted by $\mathcal{N}$), followed by a potentially empty alternating sequence from $\{\mathcal{L},\mathcal{M}\}$. 

A potentially empty alternating sequence from $\{\mathcal{L},\mathcal{M}\}$ has the form 
$$\frac{1+\mathcal{L}}{1-\mathcal{M}\mathcal{L}}+\frac{\mathcal{M}+\mathcal{L}\mathcal{M}}{1-\mathcal{L}\mathcal{M}} = \frac{(1+\mathcal{L})(1+\mathcal{M})}{1-\mathcal{L}\mathcal{M}}\,.
$$
Since we instead have a sequence of such expressions combined with a mix, we get the desired form for $\mathcal{P}_3(z)$.
\qed

\smallskip

\begin{Example}
For example, the permutation $\sigma = (6,7,4,5,9,3,8,1,10,2)$ corresponds to the word $w = \mathcal{M}\mathcal{L}\mathcal{N}\mathcal{N}$, with permutations $(1),(12),(2413),(3142)$. This permutation is shown in Figure~\ref{s:3:fig:2}.

\begin{figure}
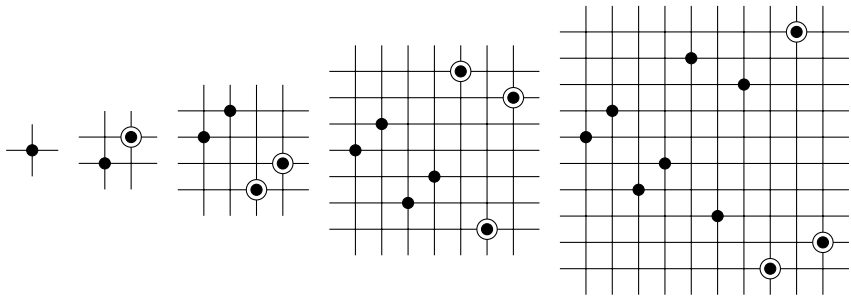

$$
\mmpattern{scale=1}{1}{1/1}{}{}
\decpatternww{scale=1}{2}{1/1,2/2}{}{2/2}{}{}{}
\decpatternww{scale=1}{4}{1/3,2/4,3/1,4/2}{}{3/1,4/2}{}{}{}
\decpatternww{scale=1}{7}{1/4,2/5,3/2,4/3,5/7,6/1,7/6}{}{5/7,6/1,7/6}{}{}{}
\decpatternww{scale=1}{10}{1/6,2/7,3/4,4/5,5/9,6/3,7/8,8/1,9/10,10/2}{}{8/1,9/10,10/2}{}{}{}
$$
\caption{A permutation in $Av(\textbf{4123,1342})$ created by adding a sum-component of (1), a skew-component of (12), a mix of (2413), and a mix of (3142).}
\label{s:3:fig:2}
\end{figure}
\end{Example}

\smallskip

Note, to label an element in $\mathcal{N}$, we need to know how the elements in the mix are arranged, and also how many of them are above $\sigma$. This means if we are adding $k$ elements, labeling them with a permutation of length $k+1$ is necessary. 

Now, it suffices to calculate $\mathcal{L},\mathcal{M}$, and $\mathcal{N}$. Clearly, $$\mathcal{L}= t-1\,,$$ where $t$ is the Catalan generating function, since new skew-components will need to avoid \textbf{123}. Note that we need $t-1$ rather than $t$, since we want nonempty elements in $Av(\textbf{123})$.

Similarly, new sum-components will need to be non-empty elements of $Av(\textbf{4123,231})$. This class was enumerated by West~\cite{West}, and the generating function $\mathcal{M}(z)$ has the expression
$$\mathcal{M}(z) = \frac{z-2z^2+2z^3}{1-4z+5z^2-3z^3}\,.$$

The mixes counted by $\mathcal{N}$ are more complicated and involve several cases. The following theorem gives our result.

\smallskip

\begin{thm}\label{s:3:thm:n}

The generating function for mixes is given by $$\mathcal{N}(z) = t^2z^3\left(t^2+\frac{z}{(1-z)^2(1-2z)}+\frac{t}{1-z}+\frac{tz}{(1-z)^3}\right)\,,$$
where $t$ is the generating function for the Catalan numbers, $$t=\frac{1-\sqrt{1-4z}}{2z}\,.$$

\end{thm}

\smallskip

Observe that $z^3$ divides $\mathcal{N}$, implying that any mix that we add must have at least three elements. This makes sense, since adding only one or two elements would be instead counted by $\mathcal{L}$, $\mathcal{M}$, or both. We want to ensure that any mix cannot be decomposed into smaller components or mixes. 

\begin{figure}
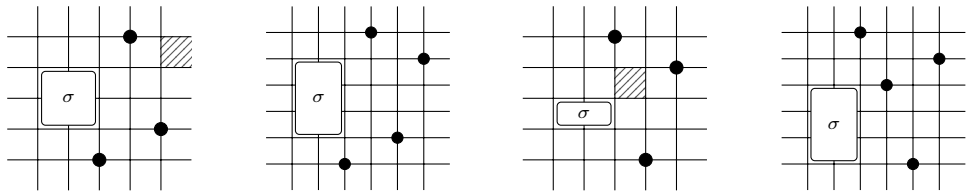

$$
\mmpattern{scale=1.17}{5}{3/1,4/5,5/2}{5/4/6/5}{1/2/3/4/\sigma}
\qquad
\mmpattern{scale=1}{6}{3/1,4/6,5/2,6/5}{}{1/2/3/5/\sigma}
\qquad
\mmpattern{scale=1.17}{5}{3/5,4/1,5/4}{3/3/4/4}{1/2/3/3/\sigma}
\qquad
\mmpattern{scale=1}{6}{3/6,4/4,5/1,6/5}{}{1/1/3/4/\sigma}
$$
\caption{Mixes counted by $\mathcal{N}_1,\mathcal{N}_2,\mathcal{N}_3,\mathcal{N}_4$, respectively. In the first and third permutations, the shaded cells must be empty.}
\label{s:3:fig:3}
\end{figure}

Before proving Theorem~\ref{s:3:thm:n}, we first describe the different ways a mix could be added to a permutation in $Av(\textbf{4123},\textbf{1342})$. When adding a mix $\tau$ to a permutation $\sigma \in Av(\textbf{4123},\textbf{1342})$, we refer to $\sigma$ as the core, and elements of the mix as $\tau_1,\tau_2,\ldots$. We count mixes based on four cases. Case 1 consists of mixes with $\tau_1$ below the core, and $\tau_1$ not serving as the 1 in a \textbf{31524}. Case 2 consists of mixes with $\tau_1$ below the core, and $\tau_1$ serving as the 1 in a \textbf{31524}. Case 3 consists of mixes with $\tau_1$ above the core, and $\tau_1$ not serving as the 5 in a \textbf{25314}. Case 4 consists of mixes with $\tau_1$ above the core, and $\tau_1$ serving as the 5 in a \textbf{25314}. We use $\mathcal{N}_1, \mathcal{N}_2, \mathcal{N}_3$, and~$\mathcal{N}_4$ as the generating functions for permutations counted by these four cases. In Figure~\ref{s:3:fig:3}, we see the minimal examples of mixes in each of the four cases.

We now describe the enumeration of each of the four cases in the following lemmas.

\smallskip

\begin{lem}\label{s:3:lem:n1}
The generating function $\mathcal{N}_1$ is equal to $$\mathcal{N}_1(z) = t^4z^3\,,$$
where $t$ is the generating function for the Catalan numbers.
\end{lem}

\smallskip

\proof
When adding a mix counted by $\mathcal{N}_1$, we must have $\tau_1$ below the core. Also, $\tau_1$ must be lower than some later element $\tau_3$ which is also below the core, since otherwise $\tau_1$ itself would be counted as a skew-component. Finally, these two elements must be split by some element above the core, $\tau_2$, since otherwise $\tau_1$ and $\tau_3$ together would form a skew-component. 

\begin{figure}
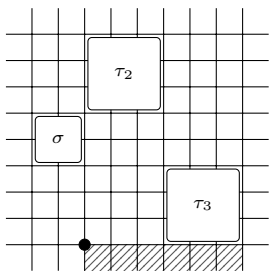

$$
\mmpattern{scale=1}{9}{3/1}{3/0,4/0,5/0,6/0,7/0,8/0}{1/4/3/6/\sigma,3/6/6/9/\tau_2,6/1/9/4/\tau_3}
$$
\caption{Structure of a mix in $\mathcal{N}_1$.}
\label{s:3:fig:n1}
\end{figure}

Observe that elements which serve as $\tau_2$ and $\tau_3$ both must avoid \textbf{12}, since otherwise a \textbf{1342} or a \textbf{4123} would occur. Other elements which are not in $\tau_2$ or $\tau_3$ could be above the core and to the right of $\tau_3$, though these belong in a mix counted by $\mathcal{N}_2$. Since we are avoiding~\textbf{4123} still, now we are in the situation where we can apply Lemma~\ref{l:tgen}. Here, every position to the right of our current minimum is available for the placement of a new minimum, except for the extreme lower right position. The reason we do not count the extreme lower right as an allowable position for a new minimum, is that a minimum there would be part of a skew-component instead of a mix. Therefore, our power of $t$ is one less than it was in our analysis of $Av(\textbf{4123},\textbf{1324})$. Because of this, our $\mathcal{R}$ for application of Lemma~\ref{l:tgen} is 
$$\mathcal{N}_1(z) = z \times \frac{tz}{1-tz} \times \frac{tz}{1-tz} = t^4z^3\,,$$
as desired.
\qed

\smallskip

Now, we consider the second case; mixes which have $\tau_1$ below the core, and serving as the 1 in an occurrence of \textbf{31524}. These mixes are counted by $\mathcal{N}_2$.

\smallskip

\begin{lem}\label{s:3:lem:n2}
The generating function $\mathcal{N}_2$ is equal to 
$$
\mathcal{N}_2(z) = \frac{t^2z^4}{(1-z)^2(1-2z)}\,,
$$
where $t$ is the generating function for the Catalan numbers.
\end{lem}

\begin{figure}
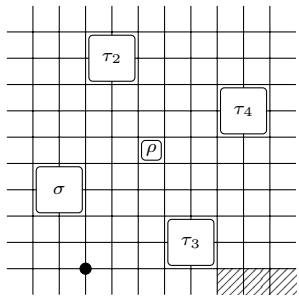

$$
\mmpattern{scale=1}{10}{3/1}{9/0,10/0,8/0}{1/3/3/5/\sigma,3/8/5/10/\tau_2,6/1/8/3/\tau_3,8/6/10/8/\tau_4,5/5/6/6/\rho}
$$
\caption{Structure of a mix in $\mathcal{N}_2$, with $\tau_2,\tau_3,\tau_4$ and $\rho$ in $Av(\textbf{12})$, with $\rho$ the only possible empty one.}
\label{s:3:fig:n2}
\end{figure}

Our initial picture is the same as $\mathcal{N}_1$, except now we have a fourth element, $\tau_4$, acting as the 4 in a \textbf{31524} with $\tau_1$ as the 1. Elements which are in the position of $\tau_4$ must again avoid~\textbf{12}, since otherwise a \textbf{4123} would be formed. There is also a possibility of having a decreasing sequence, denoted by $\rho$, between the core and $\tau_4$ vertically, and between $\tau_2$ and $\tau_3$ horizontally. The structure is shown in Figure~\ref{s:3:fig:n2}.

We want to apply Lemma~\ref{l:tgen} to count these mixes, though we need to be sure that the elements which are not minima are able to be added to the right of the $\tau_4$ block. Because $\tau_3$ is below and to the left of $\tau_4$, we cannot add a new source graph of size more than 1 unless the minimum element is to the right of $\tau_3$. 

Before we add such a source graph, we note that  with the minimum directly to the right of~$\tau_1$, and any additional elements between $\tau_3$ and $\tau_4$. We can add layers of this form indefinitely, though they do not create any additional spaces for new minima. 

Therefore, our generating function $\mathcal{R}$ before we apply Lemma~\ref{l:tgen} is 
$$
\mathcal{R}(t,z) = z \times \frac{z}{1-z} \times \frac{1}{1-z} \times \frac{z}{1-z} \times \frac{tz}{1-tz} \times \frac{1}{1-\frac{z}{1-z}}\,,
$$
with the terms coming from $\tau_1, \tau_2, \rho, \tau_3, \tau_4$, and a potentially empty sequence of source graphs with minima directly to the right of $\tau_1$. Combining, we get 
$$\mathcal{R}(t,z) = \frac{tz^4}{(1-z)^2(1-2z)(1-tz)}\,.$$
Applying Lemma~\ref{l:tgen}, we see that $\mathcal{N}_2(z) = \frac{t^2z^4}{(1-z)^2(1-2z)}$, as desired.
\qed

\smallskip

We now consider mixes with $\tau_1$ above the core, not serving as the 5 in an occurrence of~\textbf{25314}. These are counted by $\mathcal{N}_3(z)$, which is described in the following lemma.

\smallskip

\begin{lem}\label{s:3:lem:n3}
The generating function $\mathcal{N}_3(z)$ is given by 
$$
\mathcal{N}_3(z) = \frac{t^3z^3}{1-z}\,,
$$
where $t$ is the Catalan generating function.
\end{lem}

\smallskip

\begin{figure}
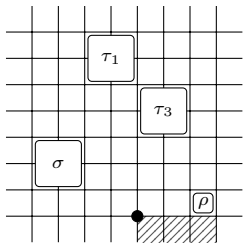

$$
\mmpattern{scale=1}{8}{5/1}{5/0,6/0,7/0}{1/2/3/4/\sigma,3/6/5/8/\tau_1,5/4/7/6/\tau_3,7/1/8/2/\rho}
$$
\caption{Structure of a mix in $\mathcal{N}_3$, with $\tau_1, \tau_3$, and $\rho$ in $Av(\textbf{12})$, with $\rho$ the only possible empty one.}
\label{s:3:fig:n3}
\end{figure}

\proof
In order for $\tau_1$ not to be part of a sum-component, there must be some element $\tau_2$ below the core to the right of $\tau_1$, and some element $\tau_3$ above the core but below $\tau_1$. This gives us the initial structure of a core counted by $\mathcal{N}_3$, as seen in Figure~\ref{s:3:fig:n3}.

Observe that the elements in $\tau_1$ form a permutation in $Av(\textbf{12})$ since otherwise a \textbf{1342} is created. Similarly, $\tau_3$ must be in $Av(\textbf{12})$, since otherwise a \textbf{4123} is created. The space to the right of $\tau_3$ and above $\tau_2$ could potentially contain a decreasing sequence, denoted by $\rho$. 

After this initial structure, we can add new source graphs, with any position to the right of~$\tau_2$ an option for new minima. Therefore, we can apply Lemma~\ref{l:tgen} with $\mathcal{R}(t,z)$ equal to 
$$
\mathcal{R}(t,z) = \frac{z}{1-z} \times z \times \frac{tz}{1-tz} \times \frac{1}{1-tz}\,,
$$
with the terms corresponding to $\tau_1,\tau_2,\tau_3$, and $\rho$, respectively. Lemma~\ref{l:tgen} then gives us the result.
\qed

\smallskip

Finally, we consider mixes counted by $\mathcal{N}_4$, with $\tau_1$ above the core, and $\tau_1$ serving as the 5 in a \textbf{25314}. The following lemma describes the enumeration of $\mathcal{N}_4$.

\smallskip

\begin{lem}\label{s:3:lem:n4}
The generating function $\mathcal{N}_4(z)$ is given by 
$$
\mathcal{N}_4(z) = \frac{t^3z^4}{(1-z)^3}\,,
$$
where $t$ is the Catalan generating function.
\end{lem}

\smallskip

\begin{figure}
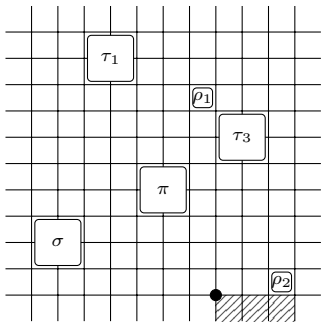

$$
\mmpattern{scale=1}{11}{8/1}{8/0,9/0,10/0}{1/2/3/4/\sigma,3/9/5/11/\tau_1,5/4/7/6/\pi,7/8/8/9/\rho_1,8/6/10/8/\tau_3,10/1/11/2/\rho_2}
$$
\caption{The structure of a mix in $\mathcal{N}_4$, where $\tau_1,\pi,\tau_3,\rho_1,\rho_2 \in Av(\textbf{12})$, with $\rho_1,\rho_2$ possibly empty.}
\label{s:3:fig:n4}
\end{figure}

\proof
We count the relevant mixes by starting with the \textbf{25314} which we know must occur. Figure~\ref{s:3:fig:n4} shows the structure of a mix counted by $\mathcal{N}_4$.

The structure is very similar to $\mathcal{N}_3$, though we now have at least one element, denoted by $\pi$, serving as the 3 in a \textbf{25314}, with $\tau_1$ as the 5. Here $\pi$ must be nonempty, and must avoid $\textbf{12}$.

Translating into generating functions, and applying Lemma~\ref{l:tgen}, we get 
$$
\aligned
\mathcal{N}_4(z) &= \frac{z}{1-z} \times \frac{z}{1-z} \times \frac{1}{1-z} \times z \times \frac{tz}{1-tz} \times \frac{1}{1-tz}\\
&= \frac{t^3z^4}{(1-z)^3}\,,
\endaligned
$$
where our terms come from $\tau_1,\pi,\rho_1,\tau_2,\tau_3$, and $\rho_2$, respectively, and $t$ is the Catalan generating function. This completes the proof.
\qed

\smallskip

We are ready to prove Theorem~\ref{s:3:thm:n}.

\smallskip

\proof[Proof of Theorem~\ref{s:3:thm:n}.] 
Since every mix comes from either $\mathcal{N}_1,\mathcal{N}_2,\mathcal{N}_3$, or $\mathcal{N}_4$, we have 
$$
\mathcal{N}(z) = \mathcal{N}_1(z)+\mathcal{N}_2(z)+\mathcal{N}_3(z)+\mathcal{N}_4(z)\,.
$$
Applying Lemmas~\ref{s:3:lem:n1}, \ref{s:3:lem:n2}, \ref{s:3:lem:n3}, and~\ref{s:3:lem:n4} and simplifying give us the desired expression for $\mathcal{N}(z)$, completing the proof of Theorem~\ref{s:3:thm:n}.
\qed

\smallskip

We are finally ready to prove Theorem~\ref{s:3:t:1}.

\smallskip

\proof[Proof of Theorem~\ref{s:3:t:1}]
Using Lemma~\ref{s:3:lem:struct}, and plugging in for $\mathcal{L},\mathcal{M}$, and $\mathcal{N}$, we find that 
$$
\mathcal{P}_3(z) = \frac{z(1-z)(1-2z)\left(\sqrt{1-4z}(1-3z+3z^2)+1-7z+17z^2-16z^3+4z^4\right)}{2-22z+96z^2-220z^3+282z^4-196z^5+64z^6-8z^7}\,,
$$
as desired.
\qed

\smallskip

Analyzing the asymptotics for $Av(\textbf{4123},\textbf{1342})$ is more difficult than our previous classes, since the denominator has an irrational root smaller than $\rho=\frac{1}{4}$, specifically at $z=0.239788$. Still, applying Proposition~\ref{prop:asymp} gives us the following result.

\smallskip

\begin{cor}\label{s:3:cor:asympP3}
The growth rate of $Av(\textbf{4123},\textbf{1342})$ is $r = \frac{1}{z} = 4.17035$. 
\end{cor}

\smallskip

The first few terms of the series expansion of $\mathcal{P}_3(z)$ are
$$1+z+2z^2+6z^3+22z^4+87z^5+352z^6+1434z^7+5861z^8+24019z^9+98677z^{10}+\ldots
$$
Subsequent terms can be found at A165533 in the OEIS~\cite{Slo}.

\bigskip

\section{Remarks}\label{s:remarks}

\medskip

\subsection{}\label{s:remarks:ss:method}
It would be interesting to analyze the applicability of enumerating permutation classes through  generating their elements by adding a sequence of source graphs. More generally, enumerating permutations by generating them through some sequence of well-defined operations which do not consist of simply adding minima one at a time. For $Av(\textbf{4123},\textbf{1342})$, generating permutations through adding sum-components, skew-components, and mixes allowed us to find generating functions. For $Av(\textbf{4123},\textbf{1243})$, decomposing permutations based on source flags rather than source graphs was fruitful. 

However, all the classes discussed in this paper avoided \textbf{4123}, meaning that adding a new left-to-right minimum required the source graph to be a fan. This is a strict condition, which definitely does not apply to most permutation classes. For classes where the basis elements do not start or end with their largest elements, the techniques we use will need to be adjusted. For example, when applying these methods to $Av(\textbf{1324})$, we see that our initial source graphs are \textbf{213}-avoiders. This is a good start, though keeping track of ways subsequent source graphs can be added appears to be much more difficult. 

\medskip

\subsection{}\label{s:remarks:ss:remaining}
The results of this paper leave only five two-by-four classes unenumerated, three of which do not appear to have algebraic generating functions (see~\cite{AHPSV}). The remaining two classes are $Av(\textbf{3412},\textbf{4123})$ and $Av(\textbf{3412},\textbf{2413})$. At first glance, the class $Av(\textbf{3412},\textbf{4123})$ is not as amenable to the methods we have used, since the initial source graph can be any permutation in the whole class, rather than being restricted to a subclass. We plan to investigate these classes further.

\medskip

\subsection{}\label{s:remarks:ss:ratio}
In~\cite{Bev2}, David Bevan posed the question of calculating when the ratio of elements in a class $Av(B)$ which avoid another pattern $\beta$ approaches a nonzero value, as $n \to \infty$. We have not answered this question, though we have found more instances of when this situation occurs. 

In Corollary~\ref{s:1:cor:ratios} we saw that as $n$ gets larger, the ratio of permutations in $Av(\textbf{4123},\textbf{1324})$ of length $n$ which also avoid \textbf{31524} is asymptotically equal to $1-\frac{45}{64n}$. 

\begin{figure}
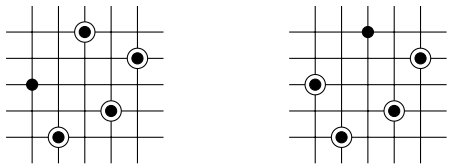

$$
\decpatternww{scale=1}{5}{1/3,2/1,3/5,4/2,5/4}{}{2/1,3/5,4/2,5/4}{}{}{}
\qquad\qquad
\decpatternww{scale=1}{5}{1/3,2/1,3/5,4/2,5/4}{}{1/3,2/1,4/2,5/4}{}{}{}
$$
\caption{The pattern \textbf{31524} contains both \textbf{1423} and \textbf{3124}.}
\label{s:4:fig:31524}
\end{figure}

The fact that this ratio tends to 1 appears to rely heavily on the structure of permutations in $Av(\textbf{4123},\textbf{1324})$. The class $Av(\textbf{4123},\textbf{1324})$ has a large subclass which avoids both \textbf{1423} and \textbf{3124}, and both these patterns are contained in \textbf{31524} (see Figure~\ref{s:4:fig:31524}). Therefore, if a permutation contains an occurrence of \textbf{31524}, then it must also contain both \textbf{1423} and \textbf{3124}, so cannot be in the large subclass. 

\vskip.5cm

\noindent
\textbf{Acknowledgments:} \. The author is grateful to 
David Bevan for presenting related results in a very understandable
way, and to Sergi Elizalde, Jay Pantone, and Henning Ulfarsson, for 
suggestions and helpful discussions. Michael Albert's PermLab software~\cite{Alb} 
was extremely useful for visualizing the structure of these permutation classes, 
and Maple~\cite{Map} was used for algebraic manipulation.

\end{document}